\documentclass{article}

\usepackage{hyperref}

\usepackage{amsfonts}
\usepackage{amssymb}
\usepackage{textcomp}
\usepackage{latexsym,amsmath}
\usepackage{graphics}

\renewcommand{\phi}{\varphi}
\newcommand{\be}{\begin{equation}}
\newcommand{\ee}{\end{equation}}
\newcommand{\ba}{\begin{eqnarray}}
\newcommand{\ea}{\end{eqnarray}}
\newcommand{\ban}{\begin{eqnarray*}}
\newcommand{\ean}{\end{eqnarray*}}

\newcommand{\nul}{{\bf0}}
\newcommand{\rd}{{\mathbb R}^d}
\newcommand{\zd}{{\mathbb Z}^d}
\newcommand{\td}{{\mathbb T}^d}

\newcommand{\z} {{\mathbb Z}}

\newcommand{\n} {{\mathbb N}}

\renewcommand{\lll}{\left(}
\newcommand{\rrr}{\right)}
\newcommand{\h}{\widehat}
\newcommand{\w}{\widetilde}
\newcommand{\too}{\mathop{\longrightarrow}}

\def\supp{\operatorname{supp}}

\def\N{{{\Bbb N}}}
\def\Z{{{\Bbb Z}}}
\def\T{{{\Bbb T}}}
\def\R{{\Bbb R}}

\def\vp{{\varphi}}
\def\t{{\theta }}

\def\({\left(}
\def\){\right)}

\def\l{{\lambda }}

\def\b{{\beta}}

\def\vp{{\varphi}}
\def\t{{\theta }}

\newtheorem{theo}{Theorem}
\newtheorem{lem}[theo]{Lemma}

\newtheorem {defi} [theo] {Definition}
\newtheorem {rem} [theo] {Remark}

\title{Approximation by  \\ sampling-type operators  in $L_p$-spaces
}
\author{
Yu. Kolomoitsev$^{1, 2}$ and M. Skopina$^{3}$
}
\date{
\small
$^{1}$Universit\"at zu L\"ubeck,
Institut f\"ur Mathematik, L\"ubeck, Germany; kolomoitsev@math.uni-luebeck.de  
\\
\small 
$^{2}$Institute of Applied Mathematics and Mechanics of NAS of Ukraine, Slov'yans'k, Ukraine
\\
$^{3}$St. Petersburg State University,  Saint Petersburg,  Russia;
skopina@ms1167.spb.edu
}

\begin{document}

\maketitle

\begin{abstract}
Approximation properties of the sampling-type quasi-projection operators  $Q_j(f,\phi, \w\phi)$ for functions $f$ from
anisotropic Besov spaces are studied. Error estimates in $L_p$-norm are obtained for a large class of tempered distributions $\w\phi$ and a large  class of functions $\phi$ under the assumptions that $\phi$ has enough decay, satisfies the Strang-Fix conditions  and a  compatibility condition with $\w\phi$. The estimates are given in terms of moduli of smoothness and best approximations.
\end{abstract}

\bigskip

\textbf{Keywords.}  Sampling-type operators, Error estimate,  Best approximation, Moduli of smoothness, Anisotropic Besov space.

\medskip

\textbf{AMS Subject Classification.} 	41A35, 41A25, 41A17, 41A15, 42B10, 94A20, 97N50

\section{Introduction}

The well-known sampling theorem (Kotel'nikov's or Shannon's formula)
states that
\be
f(x)=\sum_{k\in\z}  f(-2^{-j}k)\,\frac{\sin\pi(2^jx+k)}{\pi (2^jx+k)}
\label{0}
\ee
for  band-limited to $[-2^{j-1},2^{j-1}]$ signals (functions) $f$.
Equality~(\ref0)  holds only for functions $f\in L_2(\R)$
whose Fourier transform is supported on $[-2^{j-1},2^{j-1}]$.
However the right hand side of~(\ref0) (the sampling expansion of $f$)
has meaning for every continuous $f$ with a good enough decay, which
attracted mathematicians to study approximation properties of
classical sampling and sampling-type expansions  as $j\to\infty$. The most general form of such expansions looks as follows
$$
 Q_j(f,\phi,\w\phi)=|\det M|^{j}  \sum_{k\in\zd} \langle f,\w\vp(M^j\cdot+k)\rangle \vp(M^j\cdot+k),
$$
where $\phi$ is a function and $\w\phi$ is a distribution or function,  $M$ is a matrix, and the inner product
$\langle f,\w\vp(M^j\cdot+k) \rangle$
has meaning in some sense. The operator $Q_j(f,\phi,\w\phi)$ is usually called a quasi-projection or sampling operator.
The class of quasi-projection  operators is very
large. In particular, it includes operators with a regular function $\w\phi$,
such as scaling expansions associated with wavelet constructions
(see~\cite{BDR,  DB-DV1, v58, Jia2, KS, Sk1}  and others) and
Kantorovich-Kotelnikov operators  and their generalizations (see, e.g.,~\cite{CV0, CV2, KS3, VZ2, OT15}).  One has an essentially different class of operators $Q_j(f,\phi,\w\phi)$ if $\w\phi$ is a distribution.
This class includes the classical sampling operators, where  $\w\phi$ is the Dirac delta-function.
There are a lot of papers devoted to the study of approximation properties of such operators for different classes of functions $\phi$
(see, e.g.,~\cite{Butz4, Brown, BD, Butz5, KKS, KS, SS, Si2, Unser} and the references therein).
Consideration of  functions $\phi$ with a good decay, in particular, compactly supported ones, is very useful for apllications.
The case, where $\phi$ is a certain  linear combination of $B$-splines and $\w\phi$ is the Dirac delta-function, was studed, e.g.,
in~\cite{Butz6, Butz7, SS, Si2}. For a wide class of fast decaying functions $\phi$ and a wide class of tempered distributions $\w\phi$,
quasi-projection operators were considered in~\cite{KS}, in which error estimations in $L_p$-norm,
$2\le p\le\infty$,  were given in terms of the Fourier transform of the approximated function $f$.

The goal of the present paper is to improve the results of~\cite{KS} in several directions. In particular,
we obtain error estimates  for the case $1\le p<2$.

\section{Notation}

We use the standard multi-index notations.
    Let $\n$ be the set of positive integers, $\rd$ be the $d$-dimensional Euclidean space,
    $\zd$ be the integer lattice  in $\rd$,
    $\td=\rd\slash\zd$ be the $d$-dimensional torus.
    Let  $x = (x_1,\dots, x_d)^{T}$ and
    $y =(y_1,\dots, y_d)^{T}$ be column vectors in $\rd$,
    then $(x, y):=x_1y_1+\dots+x_dy_d$,
    $|x| := \sqrt {(x, x)}$; $\nul=(0,\dots, 0)^T\in \rd$;  		
		$\z_+^d:=\{x\in\zd:~x_k\geq~{0}, k=1,\dots, d\}.$
		If  $r>0$, then $B_r$
		denotes the ball of radius $r$ with the center in $\nul$.

  If $\alpha\in\zd_+$, $a,b\in\rd$, we set
    $$[\alpha]=\sum\limits_{j=1}^d \alpha_j, \quad
    \alpha!=\prod\limits_{j=1}^d\alpha_j!,\quad
    a^b=\prod\limits_{j=1}^d a_j^{b_j}, \,\,
    D^{\alpha}f=\frac{\partial^{[\alpha]} f}{\partial x^{\alpha}}=\frac{\partial^{[\alpha]} f}{\partial^{\alpha_1}x_1\dots
    \partial^{\alpha_d}x_d}.$$

If $A$ is a $d\times d$ matrix,
then $\|A\|$ denotes its operator norm in $\rd$; $A^*$ denotes the conjugate matrix to $A$;
 the identity matrix is denoted by $I$.

 A $d\times d$ matrix $M$ whose
eigenvalues are bigger than 1 in modulus is called a  dilation matrix.
Throughout the paper we
consider that such a matrix $M$ is fixed and  $m=|\det M|$  unless other is specified.
Since the spectrum of the operator $M^{-1}$ is
located in $B_r$, where $r=r(M^{-1}):=\lim_{j\to+\infty}\|M^{-j}\|^{1/j}$ is
the spectral radius of $M^{-1}$, and there exists at least
one point of the spectrum on the boundary of the ball, we have
	$$
	\lim_{j\to\infty}\|M^{-j}\|=0.	
	$$
A  matrix $M$ is called isotropic
	if it is similar to a diagonal matrix
such that its eigenvalues $\lambda_1,\dots,\lambda_d$ are placed on the main diagonal
	and $|\lambda_1|=\cdots=|\lambda_d|$.
	Note that if the matrix $M$ is isotropic then
	$M^*$ and $M^j$ are isotropic for all $j\in\z.$	
	It is well known that for an isotropic matrix $M$ and for  any $j\in\z$
	we have
	\be
	C^M_1 |\lambda|^j \le \|M^j\| \le C^M_2 |\lambda|^j,\quad j\in\z,
	\label{10}
\ee
   where $\lambda$ is one of the eigenvalues of $M$ and the constants $C^M_1, C^M_2$ do not depend on $j$.

By $L_p$ we  denote the space $L_p(\rd)$, $1\le p\le\infty$, with the usual norm
$\|f\|_p=\left(\int_{\rd}|f|^p\right)^{1/p}$ for $p<\infty$, and $\|f\|_\infty={\rm ess\,sup }\,|f|$.

If $f, g$ are functions defined on $\rd$ and $f\overline g\in L_1$,
then  $$\langle  f, g\rangle:=\int\limits_{\rd}f\overline g.$$

If $f\in L_1$,  then its Fourier transform is
$$
\mathcal{F}f(\xi)=\widehat
f(\xi)=\int\limits_{\rd} f(x)e^{-2\pi i
(x,\xi)}\,dx.
$$

Denote by $\mathcal{S}$ the Schwartz class of functions defined on $\rd$.
    The dual space of $\mathcal{S}$ is $\mathcal{S}'$, i.e. $\mathcal{S}'$ is
    the space of tempered distributions.
       Suppose $f\in \mathcal{S}$, $\phi \in \mathcal{S}'$, then
    $\langle \phi, f\rangle:= \overline{\langle f, \phi\rangle}:=\phi(f)$.
    If  $\phi\in \mathcal{S}',$  then $\h \phi$ denotes its  Fourier transform
    defined by $\langle \h f, \h \phi\rangle=\langle f, \phi\rangle$,
    $f\in \mathcal{S}$.
		
If $\phi$ is a function defined on $\rd$, we set
$$
\phi_{jk}(x):=m^{j/2}\phi(M^jx+k),\quad j\in\z,\quad k\in\rd.
$$
If  $\w\phi\in \mathcal{S}'$, $j\in\z, k\in\zd$, then we define $\w\phi_{jk}$ by
        $$
        \langle f, \w\phi_{jk}\rangle:=
        \langle f_{-j,-M^{-j}k},\w\phi\rangle,\quad \forall f\in \mathcal{S}.
        $$

Denote by $\mathcal{S}_N'$, $N\ge 0$, the set of  tempered distributions $\w\phi$
	whose Fourier transform $\h{\w\phi}$ is a function on $\rd$
	such that $|\h{\w\phi}(\xi)|\le C_{\w\phi} |\xi|^{N}$
	 for almost all $\xi\notin\td$
and
	 $|\h{\w\phi}(\xi)|\le C_{{\w\phi}}$
	 for almost all $\xi\in\td$.

Let $1\le p \le \infty$. Denote by ${\cal L}_p$ the set
	$$
	{\cal L}_p:=
	\left\{
	\phi\in L_p\,:\, \|\phi\|_{{\cal L}_p}:=
	\bigg\|\sum_{k\in\zd} \left|\phi(\cdot+k)\right|\bigg\|_{L_p(\td)}<\infty
	\right\}.
	$$
	It is known (see, e.g., \cite{v58}) that  ${\cal L}_p$ is a Banach space with the norm $\|\cdot\|_{{\cal L}_p}$, and
	the following properties hold:
	${\cal L}_1=L_1,$
	$\|\phi\|_p\le \|\phi\|_{{\cal L}_p}$,
	$\|\phi\|_{{\cal L}_q}\le \|\phi\|_{{\cal L}_p}$
	for $1\le q \le p \le\infty.$ Therefore, ${\cal L}_p\subset L_p$
	and ${\cal L}_p\subset {\cal L}_q$ for $1\le q \le p \le\infty.$
	If  $\phi\in L_p$ and there exist constants $C>0$
	and $\varepsilon>0$ such that 	$|\phi(x)|\le C( 1+|x|)^{-d-\varepsilon}$ for all $x\in\rd,$
	then $\phi\in {\cal L}_\infty$.

If  $A$ is a $d\times d$ matrix and  $1\le p\le \infty$, then
$$
\mathcal{B}_{A,p}:=\{g\in L_p\,:\, \supp \h g\subset A\T^d \}
$$
is a class of functions band-limited to $A\td$, and
$$
E_{A} (f)_p:=\inf\{\Vert f-g\Vert_p\,:\, g\in \mathcal{B}_{A, p}\}
$$
is the best approximation of $f\in L_p$  by band-limited functions $g$ from the class $\mathcal{B}_{A,p}$.
Additionally, we consider the following special best approximation
$$
E_{A}^*(f)_p:=\inf\{\Vert f-g\Vert_p\,:\, g\in \mathcal{B}_{A, p}\cap L_2\},
$$
which is equivalent to $E_{A} (f)_p$ for all $1\le p<\infty$ (see Lemma~\ref{lemE} below).

We will use the following anisotropic Besov spaces with respect to the matrix~$M$. We  say that
$f\in \mathbb{B}_{p,q}^s (M)$, $1\le p\le\infty$, $0<q\le \infty$, and $s>0$, if $f\in L_p$ and
$$
\Vert f\Vert_{\mathbb{B}_{p,q}^s (M)}:=\Vert f\Vert_p+\(\sum_{\nu=1}^\infty m^{\frac sd q\nu} E_{M^\nu} (f)_p^q\)^\frac1q<\infty.
$$

If $A$ is a $d\times d$  nonsingular  matrix, $f\in L_p$ and  $s\in \N$,  then
$$
\Omega_s(f,A^{-1})_p:=\sup_{|At|<1, t\in \R^d} \Vert \Delta_t^s f\Vert_p,
$$
where
$$
\Delta_t^s f(x):=\sum_{\nu=0}^s (-1)^\nu \binom{s}{\nu} f(x+t\nu).
$$
This is the so-called (total) anisotropic modulus of smoothness associated with $A$. The classical modulus of smoothness is defined
by
$$
\omega_s(f,h)_p:=\sup_{|t|<h} \Vert \Delta_t^s f\Vert_p,\quad h>0.
$$

As usual, if $\{a_k\}_{k}$ is a sequence, then
$$
\left\Vert \{a_k\}_{k}\right\Vert_{\ell_{p}}:=\left\{
                                                 \begin{array}{ll}
                                                  \displaystyle\bigg(\sum\limits_{k\in \Z^d}|a_k|^p\bigg)^\frac1p, & \hbox{if $1\le p<\infty$,} \\
                                                   \displaystyle\sup\limits_{k\in \Z^d}|a_k|, & \hbox{if $p=\infty$.}
                                                 \end{array}
                                               \right.
$$

Finally, for an appropriate  function $h\,:\, \R^d \to \mathbb{C}$ we denote the Fourier multiplier-type operator $\Lambda_h$ by
$$
\Lambda_h (f) := \mathcal{F}^{-1}(h \h f ),\quad f\in L_2.
$$

\section{ Preliminary information and auxiliary results.}
\label{sa}

Sampling-type expansions $\sum_{k\in\zd} \langle f, {\w\phi}_{jk}\rangle \phi_{jk}$ are elements of the shift-invariate
spaces generated by~$\phi$. It is well known that a function $f$ can be approximated by elements of such  spaces only if $\phi$
satisfies a special property, the so-called Strang-Fix conditions.

\begin{defi}
\label{d2}
 A  function $\phi$ is said to satisfy {\em the Strang-Fix conditions} of
order $s\in \N$  if there exists $\delta\in(0,1/2)$  such that
for any $k\in\zd,$ $k \neq \nul$, $\h\phi$ is boundedly differentiable up to order~$s$  on
	 $\{|\xi+k|<\delta\}$ and $D^{\beta}\h{\phi}(k) = 0$
	for $[\beta]<s$.
	\end{defi}

Approximation properties of the quasi-projection operators
$$
Q_j(f,\phi,\w\phi)=\sum_{k\in\zd} \langle f, {\w\phi}_{jk}\rangle \phi_{jk},
$$
with $\w\phi\in S'_N$,  were studied  in~\cite{Sk1},  \cite{KS} and~\cite{KKS}   for different classes
of  functions $\phi$. Since $\w\phi$ is a
tempered distribution, the inner product $\langle f, \w\phi_{jk}\rangle $ has meaning only for functions $f$ in $\mathcal{S}$.
To extend the class of functions $f$, the inner product
$\langle f, {\w\phi}_{jk}\rangle $ was replaced  by $\langle \h f, \h{\w\phi_{jk}}\rangle$ in these papers.
In particular, for $\phi\in {\cal L}_p$ the following result was obtained in~\cite{KS} (see Theorems 4 and 5).

	\begin{theo}
\label{theoQj}
	Let $2\le p \le \infty$, $1/p+1/q=1$, $s\in\n$,  $N\in\z_+$,  $\delta\in(0, 1/2)$, $\varepsilon>0$.  Let  also $M$ be an isotropic  dilation   matrix, $\phi \in {\cal L}_p$, and $\w\phi \in \mathcal{S}_N'$.
	Suppose
\begin{itemize}
  \setlength{\itemsep}{0cm}%
  \setlength{\parskip}{0cm}%

 	\item[$\star$]
	there exists $B_{\phi}>0$ such that
$\sum\limits_{k\in\zd}  |\h\phi(\xi+k)|^q<B_{\phi}$ for all $\xi\in\rd$;

      \item[$\star$] $\h\phi$ is boundedly differentiable up to order $s$  on
	 $\{|\xi+l|<\delta\}$ for all $l\in\zd\setminus\{\nul\}$;
the function 	$\sum\limits_{l\in\zd,\, l\neq\nul}|D^\beta \h \phi (\xi+l)|^q$
	is bounded on	$\{|\xi|<\delta\}$  for $[\beta]=s$, and
	the Strang-Fix conditions of order $s$ hold 	for $\phi$;

      \item[$\star$] $\h\phi\h{\w\phi}$ is  boundedly differentiable up to order $s$ on
	$\{|\xi|<\delta\}$, and $D^{\beta}(1-\h\phi\h{\w\phi})({\bf 0}) = 0$
	for all  $[\beta]<s$.

\end{itemize}

\noindent
 		If $f\in L_p$, $\h f\in L_q$, and
 $\h f(\xi)=\mathcal{O}(|\xi|^{-N-d-\varepsilon})$
as $|\xi|\to\infty$, then

	 $$
\bigg\|f-\sum_{k\in\zd} \langle \widehat{f}, \h{\w\phi_{jk}}\rangle \phi_{jk}\bigg\|_p\le
	 \begin{cases}
	 C |\lambda|^{-j(N+\frac dp + \varepsilon)}  &\mbox{if }
	s> N+\frac dp + \varepsilon\\
	  C (j+1)^{1/q} |\lambda|^{-js} &\mbox{if }
	 s= N+\frac dp + \varepsilon \\
	C|\lambda|^{-js}
	 &\mbox{if }
	 s< N+\frac dp + \varepsilon
	\end{cases},
$$
	where $\lambda$ is an eigenvalue of $M$ and  $C$ does not depends on $j$.
\end{theo}

There are several drawbacks in  Theorem~\ref{theoQj}. First,  an error estimate is obtained only for the case $p\ge 2$. Second, there are additional restrictions on the function $f$.
Even in the case of regular  functions $\w\phi$ for which  the inner product  $\langle{f}, {\w\phi_{jk}}\rangle $ has meaning for every  $f\in L_p$,  the error estimate is obtained only under additional assumption  $\h f(\xi)=\mathcal{O}(|\xi|^{-d-\varepsilon})$
as $|\xi|\to\infty$. Third, although Theorem~\ref{theoQj} provides approximation order for $Q_j(f,\phi,\w\phi)$,  more accurate  error estimates were not obtained in contrast to common results in approximation theory, where estimates are usually given in the terms of moduli of smoothness.

The described shortcomings of Theorem~\ref{theoQj} were avoided
in~\cite{KS3}, where the Kantorovich-type quasi-projection operators $Q_j(f,\phi,\w\phi)$ with regular functions $\w\phi$
and band-limited functions  $\phi \in {\cal B}$ were studied. 
 Here ${\cal B}$ is the class of functions $\phi$ given by
$$
\phi(x)=\int\limits_{\rd}\theta(\xi)e^{2\pi i(x,\xi)}\,d\xi,
$$
where $\theta$ is supported in a parallelepiped $\Pi:=[a_1, b_1]\times\dots\times[a_d, b_d]$ and such that	$\theta\big|_\Pi\in C^d(\Pi)$.
In particular, the following result was obtained in~\cite[Theorem~17]{KS3}.

\begin{theo}
\label{KSth2} 
Let $f\in L_p$, $1<p<\infty$, and $s\in\n$.  Suppose
 $\phi\in\cal B$, ${\rm supp\,} \h\phi\subset B_{1-\varepsilon}$ for some $\varepsilon\in (0,1)$, $\h\phi\in C^{s+d+1}(B_\delta)$ for some $\delta>0$;
 $\w\phi \in\mathcal{B}\cup \mathcal{L}_q$, $1/p+1/q=1$,  $\h{\w\phi}\in C^{s+d+1}(B_\delta)$  and  $D^{\beta}(1-\h\phi\overline{\h{\w\phi}})({\bf 0}) = 0$ for all $\beta\in\zd_+$, $[\beta]<s$.  Then
\begin{equation*}
  \bigg\|f-\sum\limits_{k\in\zd}
\langle f,\widetilde\phi_{jk}\rangle \phi_{jk}\bigg\|_p\le C \omega_s\(f,\|M^{-j}\|\)_p,
\end{equation*}
where $C$ does not depend on $f$ and $j$.
\end{theo}

A similar estimate for other classes of functions $\phi$ and $\w\phi$ follows from the results of Jia.
Namely, the next theorem can be obtained by combining Theorem~3.2 in~\cite{Jia2} with Lemma~3.2 in~\cite{Jia1}.

\begin{theo}
\label{Jia}
Let $f\in L_p$, $1\le p<\infty$ or $f\in C(\rd)$ and $p=\infty$, $M=mI$, and $s\in\n$.  Suppose
 $\phi$ and $\w\phi$ are compactly supported functions, 	the Strang-Fix conditions of order $s$ hold 	for $\phi$,  and  $D^{\beta}(1-\h\phi\overline{\h{\w\phi}})({\bf 0}) = 0$ for all $\beta\in\zd_+$, $[\beta]<s$.  Then
\begin{equation*}
  \bigg\|f-\sum\limits_{k\in\zd}
\langle f,\widetilde\phi_{jk}\rangle \phi_{jk}\bigg\|_p\le C \omega_s\(f,m^{-j}\)_p,
\end{equation*}
where $C$ does not depend on $f$ and $j$.
\end{theo}



Our goal is to avoid the mentioned above drawbacks  of  Theorem~\ref{theoQj} in the case of  $\phi\in {\cal L}_p$ and $\w\phi\in S'_N$.
For this, we need to specify the tempered distributions $\w\phi$. We will say that a  tempered distribution $\w\phi$ belongs to the class $\mathcal{S}_{N,p}'$ for some $N\ge 0$ and $1\le p\le \infty$ if $\w\phi \in \mathcal{S}_N'$ and there exists a constant $C_{\w\phi,p}$ such that
\begin{equation}\label{DefS}
  \Vert \Lambda_{{\overline{\h{\w\phi}}}}(T_\mu) \Vert_p \le C_{\w\vp,p} m^{\frac{N\mu}d} \Vert T_\mu \Vert_p
\end{equation}
for all $\mu\in \Z_+$ and $T_\mu\in \mathcal{B}_{M^\mu,p}\cap L_2$.

Obviously,~\eqref{DefS} is satisfied if $\w\phi$ is the Dirac delta-function.
If $M$ is an isotropic matrix, then  any distribution $\w\phi$ corresponding to
some differential operator (see~\cite{KS}, \cite{KKS}), namely
$$
\h{\w\phi} (\xi) = \sum_{[\beta]\le N} c_\beta (2\pi i \xi)^{\beta}, \quad N\in \z_+,\, j\in \N,
$$
is also in $\mathcal{S}_{N,p}'$. This easily follows from the well-known Bernstein inequality (see, e.g.,~\cite[p.~115]{Nik})
$$
\Vert g'\Vert_{L_p} \le \sigma \Vert g \Vert_{L_p}, \quad g\in L_p(\R),\quad \supp \h g \subset [-\sigma,\sigma].
$$


To extend the operator $Q_j(f,\phi,\w\phi)$ with $\w\phi \in \mathcal{S}_{N,p}'$
to the space $ \mathbb{B}_{p,1}^{d/p+N}(M)$, we  need the following lemmas.

\begin{lem}\label{lemMZ} {\sc (\cite[Theorem 4.3.1]{TB})}
Let $1\le p<\infty$, $f\in L_p$, and $\supp \h f\subset [-\sigma,\sigma]^d$, $\sigma>0$. Then
$$
\sum_{k\in \Z^d} \max_{x\in Q_{k,\sigma}} |f(x)|^p \le C(p,d) \sigma^d\Vert f\Vert_p^p,
$$
where $Q_{k,\sigma}=[-\frac1\sigma,\frac1\sigma]^d+\frac{2k}\sigma$. 
\end{lem}

\begin{lem} {\sc (\cite[Theorem 2.1]{v58})}
\label{lemKK1}
Let  $1\le p\le\infty$,   $\phi\in {\cal L}_p$, and  $\{a_k\}_{k\in\zd}\subset \mathbb{C}$.  Then
\begin{equation*}
  \bigg\Vert \sum_{k\in \Z^d} a_k \phi_{0k}\bigg\Vert_p
\le  \|\phi\|_{\mathcal{L}_p} \left\Vert \{a_k\}_{k}\right\Vert_{\ell_{p}}.
\end{equation*}
\end{lem}

\begin{lem} {\sc (\cite[Proposition 5]{NU})}
\label{lemNU}
Let  $1\le p\le\infty$,  $1/p+1/q=1$, $f\in L_p$, and  $\phi\in {\cal L}_q$. Then
\begin{equation*}
\left\Vert \{\langle f, \phi_{0k}\rangle\}_{k}\right\Vert_{\ell_{p}}
\le  \|\phi\|_{\mathcal{L}_q}\|f\|_p .
\end{equation*}
\end{lem}

\begin{lem}\label{lemE}
 Let $1\le p<\infty$ and $A$ be $d\times d$ matrix. Then
\begin{equation}\label{ee}
  E_{A}(f)_p\le E_{A}^*(f)_p\le C(p,d) E_{A}(f)_p.
\end{equation}
\end{lem}

{\bf Proof.}
We need to prove only the second inequality in~\eqref{ee}.
Let $Q\in \mathcal{B}_{A,p}$ be such that $\Vert f-Q\Vert_p\le 2E_A(f)_p$.
If $1\le p\le 2$, then using the embedding $\mathcal{B}_{A,p}\subset \mathcal{B}_{A,2}$ (see, e.g.,~\cite[Theorem~3.3.5]{Nik}), we get
$
E_A^*(f)_p\le \Vert f-Q\Vert_p\le 2E_A(f)_p.
$

Now consider the case $p >2$. Let $T\in \mathcal{S}$ be such that $\Vert f-T\Vert_p\le E_A(f)_p$.
Denote $S_A g:=\mathcal{F}^{-1} \chi_{A\T^d} * g$. 
It is clear that $S_A T\in L_p\cap L_2$ and $S_A Q=Q$. Thus, taking into account that $S_A$ is bounded in $L_p$ for all $1<p<\infty$ (see, e.g.,~\cite[Section~1.5.7]{Nik}), we derive
\begin{equation*}
  \begin{split}
      E_A^*(f)_p&\le \Vert f-S_A T\Vert_p\le \Vert f-Q\Vert_p+\Vert S_A(Q-T)\Vert_p\\
      &\le \Vert f-Q\Vert_p+C_1 \Vert Q-T\Vert_p\le (1+C_1)\Vert f-Q\Vert_p+C_1 \Vert f-T\Vert_p\le C_2 E_A(f)_p,
   \end{split}
\end{equation*}
which proves the lemma.~~$\Diamond$

In particular Lemma~\ref{lemE} implies that for any $f\in L_p$, $1\le p<\infty$, there exists a function $T \in \mathcal{B}_{A,p}\cap L_2$ such that $\Vert f-T_\mu \Vert_p\le C(p,d) E_{A}(f)_p$.

\begin{lem}
\label{lem1}
Let  $1\le p<\infty$, $N\ge 0$, $\delta \in (0,1]$, $\w\vp\in \mathcal{S}_{N,p}'$,
and $f\in \mathbb{B}_{p,1}^{d/p+N}(M)$. Suppose that  functions
$T_\mu \in \mathcal{B}_{\delta M^\mu,p}\cap L_2$, $\mu\in\z_+$,  are such that
$$
\Vert f-T_\mu \Vert_p\le C(p,d) E_{\delta M^\mu}(f)_p.
$$
Then the sequence $\{\{\langle \h T_\mu, \h{\w\phi}_{0k}\rangle\}_k\}_{\mu=0}^\infty$
converges in $\ell_p$ as $\mu\to\infty$;   a fortiori  the sequence
 $\{\langle \h T_\mu, \h{\w\phi}_{0k}\rangle\}_{\mu=0}^\infty$ converges uniformly with respect to $k\in\zd$,  and the limit does not depend on
the choice of functions~$T_\mu$;
$$
	\sum_{\mu=n}^\infty
	\Vert \{\langle \h T_{\mu+1} -\h T_\mu, \h{\w\phi}_{0k}\rangle \}_{k} \Vert_{\ell_p}\le
	 C\sum_{\mu=n}^\infty m^{\mu(\frac{N}d+\frac1p)}E_{\delta M^\mu}(f)_p,\quad \forall n\in\z_+
	$$
where the constant $C$ depends only on $\w\vp$, $d$, $p$, and $M$.
\end{lem}

{\bf Proof.} Set
$$
F(x):=\int\limits_{\rd}\Big(\h T_{\mu+1}({M^*}^{\mu+1}\xi)-\h T_{\mu}({M^*}^{\mu+1}\xi)\Big)
\overline{\h{\w \phi}({M^*}^{\mu+1}\xi)}e^{2\pi i(\xi,x)}\,d\xi.
$$
We have
\begin{equation}
\label{121}
\begin{split}
     \Vert \{\langle \h T_{\mu+1}, \h{\w\phi}_{0k}\rangle -
		\langle \h T_\mu, \h{\w\phi}_{0k}\rangle \}_{k} \Vert_{\ell_p}^p&=\sum_{k\in\zd}\bigg|\int\limits_{\rd}\Big(\h T_{\mu+1}(\xi)-\h T_{\mu}(\xi)\Big)
\overline{\h{\w \phi}(\xi)}e^{-2\pi i(\xi,k)}\,d\xi\bigg|^p\\
&=m^{p(\mu+1)}\sum_{k\in\zd}|F(M^{\mu+1}k)|^p.
\end{split}
\end{equation}	
Since
$$
\supp\h T_{\mu+1}({M^*}^{\mu+1}\cdot)\subset \td,\quad
\supp\h T_{\mu}({M^*}^{\mu+1}\cdot)\subset \sqrt d\|M^{-1}\|\td,
$$
one can  easily see that
$$
\supp \h F\subset [-\sigma, \sigma]^d,\quad \text{where}\ \  \sigma=\sigma(M,d).
$$
Thus, by Lemma~\ref{lemMZ}, taking into account that  each set
	$Q_{r,\sigma}$ contains only finite number (depending only on $M$ and $d$) of points $M^{\mu+1}k\in \zd$, $k\in\zd$, we obtain
	$$
	\sum_{k\in\zd}|F(M^{\mu+1}k)|^p\le C_1\sigma^d\int\limits_{\rd}|F(x)|^p\,dx.
	$$
After the change of variables, the above inequality yields
	$$
	m^{p(\mu+1)}\sum_{k\in\zd}|F(M^{\mu+1}k)|^p\le C_1\sigma^dm^{\mu+1}\Vert \Lambda_{\mathcal{F}^{-1}{{\w\phi}}}(\overline{T_{\mu+1}}-\overline{T_\mu}) \Vert_p^{p}.
	$$
Using this estimate together with~\eqref{121} and~\eqref{DefS}, we get
	$$
	\Vert \{\langle \h T_{\mu+1}, \h{\w\phi}_{0k}\rangle -
	\langle \h T_\mu, \h{\w\phi}_{0k}\rangle \}_{k} \Vert_{\ell_p}\le
	C_2  m^{\mu(\frac{N}d+\frac1p)}\|T_{\mu+1}-T_\mu\|_p.
	$$

Hence, for any $n\in\z_+$, we have
\begin{equation*}
  \begin{split}
      	\sum_{\mu=n}^\infty
	\Vert \{\langle \h T_{\mu+1}, \h{\w\phi}_{0k}\rangle -\langle \h T_\mu, \h{\w\phi}_{0k}\rangle \}_{k} \Vert_{\ell_p}&\le
	C_4\sum_{\mu=n}^\infty m^{\mu(\frac{N}d+\frac1p)}\|T_{\mu+1}-T_\mu\|_p\\
 &\le 2 C_4\sum_{\mu=n}^\infty m^{\mu(\frac{N}d+\frac1p)}E_{\delta M^\mu}(f)_p,
   \end{split}
\end{equation*}
and it is easy to see that $C_4$ depends only on $\w\phi$, $d$, $p$, and $M$.

The latter series is convergent because $f\in \mathbb{B}_{p,1}^{d/p+N}(M)$,
which yields that $\{\{\langle \h T_\mu, \h{\w\phi}_{0k}\rangle\}_k\}_{\mu=0}^\infty$ is a Cauchy sequence in $\ell_p$. In particular, this implies that
the sequence $\{\langle \h T_\mu, \h{\w\phi}_{0k}\rangle\}_{\mu=0}^\infty$ is convergent for every $k\in\zd$.
It remains to  check that the limit does not depend on
the choice of functions $T_\mu$.  Let $\{T'_\mu\}_\mu$ be another sequence of functions with the same properties.  Repeating all above arguments with $T'_\mu$ instead of $T_{\mu+1}$, we obtain
	$$
		\Vert \{\langle \h T'_{\mu}, \h{\w\phi}_{0k}\rangle -\langle\h T_\mu, \h{\w\phi}_{0k}\rangle \}_{k} \Vert_{\ell_p}\le
	C_5 m^{\mu(\frac{N}d+\frac1p)}\|T'_{\mu}-T_\mu\|_p\le 2 C_5 m^{\mu(\frac{N}d+\frac1p)}E_{\delta M^\mu}(f)_p,
	$$
	where $C_5$ depends only on $\w\phi$, $d$,  $p$,  and $M$.
	
It follows that $\langle \h T'_{\mu}, \h{\w\phi}_{0k}\rangle -\langle \h T_\mu, \h{\w\phi}_{0k}\rangle\to 0$ as $\mu\to\infty$ for any $k\in\zd$,
which completes the proof.~~$\Diamond$

\bigskip

Consider the functional $\langle f,\w\vp_{jk}\rangle$.   If $\w\phi\in \mathcal{S}'_N$, then  $\langle f,\w\vp_{jk}\rangle$  has meaning whenever $f\in \mathcal{S}$.
Due to Lemma~\ref{lem1} we can extend  the linear functional  $\langle f,\w\vp_{jk}\rangle$ from $\mathcal{S}$ to  $\mathbb{B}_{p,1}^{d/p+N}(M)$  as follows.

	\begin{defi}
\label{def0}
	Let $1\le p<\infty$, $N\ge 0$, $\w\phi\in S'_{N, p}$, and the functions $T_\mu$ be as in Lemma~\ref{lem1}.
	For every $f\in \mathbb{B}_{p,1}^{d/p+N}(M)$ and $k\in\zd$, we set
$$
\langle f,\w\vp_{0k}\rangle:=\lim\limits_{\mu\to\infty}\langle\h T_\mu,\h{\w\phi}_{0k}\rangle
$$
and
	$$
	\langle f,\w\vp_{jk}\rangle:= m^{-j/2}\langle f(M^{-j}\cdot),\w\vp_{0k}\rangle, \quad j\in\z_+.
	$$
	\end{defi}	
	
	Now, the quasi projection  operators
	$$
	Q_j(f,\phi,\w\phi)=\sum_{k\in\zd}\langle f,\w\vp_{jk}\rangle\phi_{jk}
	$$
	are  defined on the space $\mathbb{B}_{p,1}^{d/p+N}(M)$ for a wide class of appropriate functions $\phi$.

\begin{rem}
\label{prop002}
If $\phi\in {\cal L}_p$, then  the operator $Q_j(f,\phi,\w\phi)$ is well defined
because the series $\sum\limits_{k\in\zd}\langle f,\w\vp_{jk}\rangle\phi_{jk}$
converges unconditionally in $L_p$, which follows from Lemmas~\ref{lemKK1} and~\ref{lem1}.
\end{rem}

\begin{rem}
\label{prop003}
If the Fourier transform of a function $f$ has enough
decay for the inner product $\langle \h f,\h{{\w \vp}_{0k}}\rangle$ to have meaning, then
it is natural to define $Q_j(f,\phi,\w\phi)$ setting $\langle f,\w\vp_{0k}\rangle:=\langle \h f,\widehat{{\w \vp}_{0k}}\rangle$
(see Theorem~\ref{theoQj}).
Such an operator $Q_j(f,\phi,\w\phi)$ will be the same as one  in correspondence with Definition~\ref{def0}.
\end{rem}

{\bf Proof.}
We assume that $f\in L_p$, $\h f\in L_q$, $1/p+1/q=1$,
 $\h f(\xi)=\mathcal{O}(|\xi|^{-N-d-\varepsilon})$ as $|\xi|\to\infty$, $\varepsilon>0$,
 $f\in \mathbb{B}_{p,1}^{d/p+N}(M)$, $\w\phi\in S'_N$, and $M$ is an isotropic matrix.

First we consider the case $p<2$.
 Let $T_\mu \in\mathcal B_{\delta M^\mu, p}$
be such that
$
\|f-T_\mu\|_p\le 2E_{\delta M^\mu}(f)_p
$
Since $\mathcal{B}_{{\delta M^\mu},p}\subset \mathcal{B}_{{\delta M^\mu},2}$ (see, e.g.,~\cite[Theorem~3.3.5]{Nik}), we have $T_\mu\in\mathcal B_{\delta M^\mu, p}\cap L_2$ and, due to the Hausdorff-Young inequality,
\be
\label{301}
\|\h f-\h T_{\mu}\|_q\le \|f-T_{\mu}\|_p\le 2E_{\delta M^\mu}(f)_p,
\ee
Next we have
\begin{equation*}
  \begin{split}
      |\langle \h f,\h{\w\vp_{0k}}\rangle-\langle \h T_\mu,\h{\w\vp_{0k}}\rangle|\le C_1 &\int\limits_{|\xi|\le \delta\sqrt d\|M^\mu\|} (1+|\xi|)^N |\h f(\xi)-\h T_\mu(\xi)| d\xi
\\
&\qquad+ C_2 \int\limits_{|\xi|\ge \delta\sqrt d\|M^\mu\|} (1+|\xi|)^N |\h f(\xi)-\h T_\mu(\xi)| d\xi
=:C_1I_1+C_2 I_2.
   \end{split}
\end{equation*}
If $|\xi|\ge \delta\sqrt d\|M^\mu\|$, then $\xi\not\in \delta M^\mu\td$, and hence
$\h T_\mu(\xi)=0$, which yields
$$
I_2=\int\limits_{|\xi|\ge \delta\sqrt d\|M^\mu\|} (1+|\xi|)^N |\h f(\xi)| d\xi \too_{\mu\to\infty}0.
$$
Using H\"older's inequality and~\eqref{301}, we derive
$$
I_1\le 2\Bigg(\int\limits_{|\xi|\le \delta\sqrt d\|M^\mu\|} (1+|\xi|)^{pN}\Bigg)^{1/p}E_{\delta M^\mu}(f)_p
\le C_3\|M^\mu\|^{N+\frac dp}E_{\delta M^\mu}(f)_p.
$$
Taking into account that $f\in \mathbb{B}_{p,1}^{d/p+N}(M)$ and the inequality $\|M^\mu\|\le C(d, M) m^{\mu/d}$, which holds for any isotropic matrix,
we obtain that $I_1\to0$ as $\mu\to\infty$. Thus
\be
\label{302}
 |\langle \h f,\h{\w\vp_{0k}}\rangle-\langle \h T_\mu,\h{\w\vp_{0k}}\rangle|
\too_{\mu\to\infty}0.
\ee

Let now $p\ge2$. Denote $T_\mu=\mathcal{F}^{-1}(\chi_{\delta M^\mu\td})*{g_\mu}$, where
 $g_\mu \in \mathcal{S}$ is such that
\be
\label{303}
\Vert \h f-\h{g_\mu}\Vert_q\le E_{\delta M^\mu}(f)_p.
\ee
Obviously, $T_\mu\in\mathcal B_{\delta M^\mu, p}\cap L_2$ and
$\h T_\mu =\chi_{\delta M^\mu\td}\h{g_\mu}$.
Next we have
\begin{equation*}
  \begin{split}
       &|\langle \h f,\h{\w\vp_{0k}}\rangle-\langle \h T_\mu,\h{\w\vp_{0k}}\rangle|\\
&\le \int\limits_{\rd} (1+|\xi|)^N |\h f(\xi)-\chi_{\delta M^\mu\td}(\xi)\h f(\xi)| d\xi
+\int\limits_{\rd} (1+|\xi|)^N |\chi_{\delta M^\mu\td}(\xi)(\h f(\xi)-\h{g_\mu}(\xi)| d\xi\\
&= \int\limits_{\xi\not \in\delta M^\mu\td} (1+|\xi|)^N |\h f(\xi)| d\xi
+\int\limits_{\xi\in\delta M^\mu\td} (1+|\xi|)^N
|\h f(\xi)-\h{g_\mu}(\xi)| d\xi
   \end{split}
\end{equation*}
Due to the same arguments as above, using~\eqref{303} instead of~\eqref{301},
we conclude  again that relation~\eqref{302} holds true.
To complete the proof it remains to note that the limit of the sequence  $\{\langle \h T_\mu,\h{\w\vp_{0k}}\rangle\}_\mu$ does not depend on the choice of appropriate functions $T_\mu$.~~$\Diamond$

\medskip

\medskip

\begin{lem}
\label{thKS}
Let $1\le p< \infty$,  $1/p+1/q=1$, $\delta\in(0,1]$,  $\phi\in {\cal L}_p$, and the functions $T_\nu$, $\nu\in \Z_+$, be as in Lemma~\ref{lem1}.
%

\noindent $(i)$ If $N\ge 0$, $\w\phi\in S'_{N,p}$,  $f\in \mathbb{B}_{p,1}^{d/p+N}(M)$, and  $\nu\in\z_+$, then
\begin{equation}\label{KS000}
\begin{split}
\|f-Q_0(f,\vp,\w\vp)\|_p\le
\|T_\nu-Q_0(T_\nu,\vp,\w\vp)\|_p+
C\sum_{\mu=\nu}^\infty m^{\mu(\frac{N}d+\frac1p)}E_{\delta M^\mu}(f)_p,
\end{split}
\end{equation}
where  the constant $C$ does not depend on $f$ and $\nu$.

\noindent $(ii)$ If  $\w\vp \in \mathcal{L}_q$,
 $f\in L_p$, and  $\nu\in\z_+$, then
\begin{equation}\label{KS000_K}
\begin{split}
\|f-Q_0(f,\vp,\w\vp)\|_p\le
\|T_\nu-Q_0(T_\nu,\vp,\w\vp)\|_p+
CE_{\delta M^\nu}(f)_p,
\end{split}
\end{equation}
where  the constant $C$ does not depend on $f$ and $\nu$.
\end{lem}

{\bf Proof.}
We have
\begin{equation}\label{_1}
  \begin{split}
     \|f-Q_0(f,\vp,\w\vp)\|_p\le \|f-T_\nu\|_p+\|T_\nu-Q_0(T_\nu,\vp,\w\vp)\|_p+\|Q_0(f-T_\nu,\vp,\w\vp)\|_p.
  \end{split}
\end{equation}
Under assumptions of item $(i)$, by  Definition~\ref{def0}, using  Lemmas~\ref{lemKK1} and~\ref{lem1},  we derive
\begin{equation}\label{_2}
  \begin{split}
     \|Q_0(f-T_\nu,\vp,\w\vp)\|_p&\le C_1\|\{\langle f-T_\nu, \w\phi_{0k}\rangle\}_k\|_{\ell_p}\le C_2\sum_{\mu=\nu}^\infty\|\{\langle T_{\mu+1}-T_\mu, \w\phi_{0k}\rangle\}_k\|_{\ell_p}\\
&\le C_3\sum_{\mu=\nu}^\infty m^{\mu(\frac{N}d+\frac1p)}E_{\delta M^\mu}(f)_p.
  \end{split}
\end{equation}
Combining~\eqref{_1} and~\eqref{_2}, we get~\eqref{KS000}.

Under assumptions of item $(ii)$, it follows from Lemmas~\ref{lemKK1} and~\ref{lemNU} that
\begin{equation}\label{_3}
  \begin{split}
     \|Q_0(f-T_\nu,\vp,\w\vp)\|_p&\le \|\phi\|_{\mathcal{L}_p}\|\{\langle f-T_\nu, \w\phi_{0k}\rangle\}_k\|_{\ell_p}\\
&\le \|\phi\|_{\mathcal{L}_p}\Vert \w\vp \Vert_{\mathcal{L}_q}\Vert f-T_\nu\Vert_p\le C_4E_{\delta M^\nu}(f)_p.
  \end{split}
\end{equation}
Combining~\eqref{_1} and~\eqref{_3}, we obtain~\eqref{KS000_K}.~~$\Diamond$

\bigskip

We will use the following  Jackson-type theorem (see, e.g.,~\cite[Theorem 5.2.1 (7)]{Nik} or~\cite[5.3.2]{Timan}).

\begin{lem}\label{lemJ}
  Let $f\in L_p$, $1\le p\le \infty$, 
		and $s\in \N$. Then
\begin{equation*}
  E_{A}(f)_p\le C \Omega_{s}(f,A^{-1})_p,
\end{equation*}
where $C$ is a constant independent of $f$ and $A$.
\end{lem}

\begin{lem}\label{lem1'} {\sc (See~\cite{Wil}.)}
  Let $1\le p\le\infty$, $s\in \N$, and $T\in \mathcal{B}_{I,p}$. Then
\begin{equation*}
   \sum_{[\beta]=s} \Vert D^\beta T\Vert_p\le C\omega_s(T,1)_p,
\end{equation*}
where the constant $C$ depends only on $p$ and $s$.
\end{lem}


\section{Main results}

\begin{theo}
\label{corMOD1'+}
Let $1\le p<\infty$,  $N\ge 0$, $s\in \N$, $\delta\in (0,1/2)$.
Suppose that  $\vp \in \mathcal{L}_p$ and $\w\phi\in S'_{N,p}$   satisfy the following conditions:
\begin{itemize}
  \item[1)] 	$\h{\w\phi}\in  C^{s+d+1}({2\delta\td)}$;
  \item[2)] $\h\phi\in C^{s+d+1}\(\cup_{l\in\zd}(2\delta\td+l)\)$;
  \item[3)] the Strang-Fix conditions of order $s$  hold for $\phi$;
  \item[4)] $D^{\beta}(1-\h\phi\overline{\h{\w\phi}})({\bf 0}) = 0$ for all  $\beta\in\zd_+$, $[\beta]<s$;
  \item[5)]$ \sum\limits_{l\ne\nul}\sup\limits_{\xi\in{2\delta\td}}|D^{\beta}\h\phi(\xi+l)|\le B_\phi$ for all  $\beta\in\zd_+$, $s\le [\beta]\le s+d+1$.
\end{itemize}
Then, for every $f\in \mathbb{B}_{p,1}^{d/p+N}(M)$ and $j\in \z_+$, we have
\begin{equation}\label{110}
\begin{split}
\Vert f - Q_j(f,\phi,\w\phi) \Vert_p\le C \(\Omega_s(f, M^{-j})_p+m^{-j(\frac1p+\frac Nd)}\sum_{\nu=j}^\infty m^{(\frac1p+\frac Nd)\nu} E_{\delta M^\nu}(f)_p\),
\end{split}
\end{equation}
and if, moreover, $\w\vp \in \mathcal{L}_q$, $1/p+1/q=1$,    then for every $f\in L_p$
\begin{equation}\label{110K}
  \Vert f - Q_j(f,\phi,\w\phi) \Vert_p\le C \Omega_s(f,M^{-j})_p,
\end{equation}
where the constant $C$ does not depend on $f$ and $j$.
\end{theo}

{\bf Proof.} First of all note that it suffices to prove~\eqref{110} and~\eqref{110K} for  $j=0$. Indeed,
	$$
\bigg\|f-\sum\limits_{k\in\zd}\langle f,\widetilde\phi_{jk}\rangle \phi_{jk}\bigg\|_p=
\bigg\|m^{-j/p}f(M^{-j}\cdot)-\sum\limits_{k\in\zd}\langle m^{-j/p}f(M^{-j}\cdot),\widetilde\phi_{0k}\rangle \phi_{0k}\bigg\|_p.
$$
Obviously, $m^{-j/p}f(M^{-j}\cdot)\in \mathbb{B}_{p,1}^{d/p+N}(M^0)$ 	whenever $f\in \mathbb{B}_{p,1}^{d/p+N}(M)$. We have also  that
$$
E_{\delta M^\nu}(m^{-j/p}f(M^{-j}\cdot))_p\le E_{\delta M^{\nu+j}}(f)_p
$$
and
$$
\omega_s(f(M^{-j}\cdot), 1)_p= m^{j/p}\Omega_s(f, M^{-j})_p,
$$
 which yield~\eqref{110} and~\eqref{110K} whenever these relations hold true for $j=0$.

Let $T \in \mathcal{B}_{\delta I,p}\cap L_2$ be such that $\Vert f-T\Vert_p\le 2 E_{\delta I}(f)_p$.  Due to Lemma~\ref{thKS}, we need only to check that
\be
\label{104}
\bigg\|T-\sum\limits_{k\in\zd}\langle T,\widetilde\phi_{0k}\rangle \phi_{0k}\bigg\|_p\le
C_1\omega_s(f, 1)_p.
\ee

Since $\supp \h T\subset \delta\td$, without loss of generality, in what follows we
assume that $\h{\w\phi}\in C^{s+d+1}(\rd)$ and  $\supp\h{\w\phi}\subset 2\delta\td$.

Using Lemma~1 from~\cite{KS}, Carleson's theorem (with cubic convergence of the Fourier series),
we derive
\be
\label{102}
\mathcal{F}\Big(\sum_{k\in\zd}\langle  T, \w\phi_{0k}\rangle\phi_{0k}\Big)=G\h\phi,
\ee
where $G(\xi):=\sum_{k\in\zd}\h T(\xi+k)\overline{\h{\w\phi}(\xi+k)}$.

Consider auxiliary functions $\phi^*$ and $\phi^{(\beta)}$, $[\beta]=s$, $\b\in \Z_+^d$, such that   $\h{\phi^*}, \h{\phi^{(\beta)}}\in C^{s+d+1}(\rd)$,
$$
 \h{\phi^*}(\xi)=\left\{
\begin{array}{ll}
\h\phi(\xi), & \xi\in\delta\td,
\\
0, & \xi\in\td\setminus 2\delta\td,
 \end{array}
                \right.
$$
$$
\h{\phi^{(\beta)}}(\xi)=0\quad\text{for}\quad \xi\in\td,
$$
and for any $l\in \Z^d\setminus \{0\}$
\be
\label{107}
\h{\phi^{(\beta)}}(\xi+l)= \left\{
                  \begin{array}{lll}
             \displaystyle      \frac{s}{\beta!}  \xi^{\beta}
						\int\limits_0^1 (1-t)^{s-1} D^{\beta}\h{\phi}( t\xi+l) d t, & \xi\in\delta\td,
\\
              \displaystyle      0, & \xi\in\td\setminus 2\delta\td.
                  \end{array}
                \right.
\ee
In view of condition 5), we can also provide that
\be
\label{106}
\sum\limits_{l\ne\nul}\sup\limits_{\xi\in2\delta\td}\left|
D^{\beta'}\left(\frac{\h{\phi^{(\beta)}}(\xi+l)}{\xi^\beta}\right)\right|\le B'_\phi,
  \quad\beta'\in\zd_+, \ 0\le[\beta']\le d+1.
\ee

By condition 3), taking into account the fact that $\supp \h{T}\subset \delta \T^d$ and that, by Taylor's formula with the integral reminder,
$$
  \h{\phi} (\xi+l) = \sum_{[\beta]=s} \frac{s}{\beta!}  \xi^{\beta}
	\int\limits_0^1 (1-t)^{s-1} D^{\beta}\h{\phi}( t\xi+l) d t,
$$
we have
\be
\label{103}
G(\xi)\h\phi(\xi)=\h T(\xi)\overline{\h{\w\phi}(\xi)}\h{\phi^*}(\xi)+G(\xi)
\sum_{[\beta]=s}\h{\phi^{(\beta)}}(\xi).
\ee

It is easy to see that $D^{\beta'}\h{\phi^*},D^{\beta'}\h{\phi^{(\beta)}}\in L$ for all
$\beta'\in\zd_+$, \ $0\le[\beta']\le d+1$. It follows that  the functions  ${\phi^*}=\mathcal{F}^{-1}(\h{\phi^*})$ and ${\phi^{(\beta)}}=\mathcal{F}^{-1}(\h{\phi^{(\beta)}})$ are in $\mathcal{L}_p$.

 Using~\eqref{102} and \eqref{103}, we have
\be
\begin{split}
   \bigg\|T-\sum\limits_{k\in\zd}\langle T,\widetilde\phi_{0k}\rangle \phi_{0k}\bigg\|_p&\le\bigg\|T-\sum\limits_{k\in\zd}\langle T,\widetilde\phi_{0k}\rangle \phi^*_{0k}\bigg\|_p+\sum_{[\beta]=s}
\bigg\|\sum\limits_{k\in\zd}\langle T,\widetilde\phi_{0k}\rangle \phi^{(\beta)}_{0k}\bigg\|_p\\
&=: I^*+ \sum_{[\beta]=s}I^{(\beta)}.
\end{split}
\label{101}
\ee

By condition 4),  using the Jensen inequality and taking into account the fact that $\supp \h{T}\subset \delta \T^d$ and that by Taylor's formula with the integral reminder
\begin{equation*}
  \h{\phi} (\xi)\overline{\h{\w\phi}(\xi)} = 1 + \sum_{[\beta]=s} \frac{s}{\beta!}  \xi^{\beta} \int\limits_0^1 (1-t)^{s-1} D^{\beta}\h{\phi}\overline{\h{\w\phi}}( t \xi) d t,
\end{equation*}
we obtain
\begin{equation*}
  \begin{split}
       I^*&=\left\Vert \mathcal{F}^{-1}\left\{ \({1-(\overline{\h{\w\vp}}\h{\vp^*})}\)  \h{T}\right\}\right\Vert_{p}\\
&\le \sum_{[\beta]=s} \frac{s}{\beta!}  \int\limits_0^1\,dt
\bigg( \int\limits_{\rd}\bigg| \int\limits_{\rd} D^{\beta}\h{\phi^*}\overline{\h{\w\phi}}( t \xi)
\xi^\beta\h{T}(\xi) e^{2\pi i(x,\xi)}\,d\xi\bigg|^p\,dx\bigg)^{1/p}
\\
&=\sum_{[\beta]=s} \frac{s}{\beta!}  \int\limits_0^1\,dt
\bigg( \int\limits_{\rd}\bigg| \int\limits_{\rd}t^{-d} \mathcal{F}^{-1}(D^{\beta}\h{\phi^*}\overline{\h{\w\phi}})\Big( \frac yt\Big)
{D^\beta T}(x-y) \,dy\bigg|^p\,dx\bigg)^{1/p}.
   \end{split}
\end{equation*}
Since, obviously, the function $ \mathcal{F}^{-1}(D^{\beta}\h{\phi^*}\overline{\h{\w\phi}})$ is summable on $\rd$, this yields
\be
\label{105}
 I^*\le\sum_{[\beta]=s} \frac{s}{\beta!}\|\mathcal{F}^{-1}(D^{\beta}\h{\phi^*}
\overline{\h{\w\phi}})\|_1 \|D^\beta T\|_p
\le C_2\sum_{[\beta]=s}\|D^\beta T\|_p.
\ee

Next, we fix $\beta$ and estimate $I^{(\beta)}$. Set
$$
\h\psi(\xi)=\frac{\h{\phi^{(\beta)}}(\xi)}{\gamma(\xi)},
$$
where   $\gamma$ is a $1$-periodic functions (in each variable) such that
$\gamma\in C^{d+1}(\rd)$ and  $\gamma(\xi)=\xi^\beta$ on ${2\delta\td}$.

Using~\eqref{107} and~\eqref{106}, we have
$$
\int\limits_{\rd}|\h\psi(\xi)| \,d\xi =\sum_{k\in\zd}\int\limits_{\td-k}
\frac{|\h{\phi^{(\beta)}}(\xi)|}{|\gamma(\xi)|}d\xi
=\int\limits_{2\delta\td}\frac{\sum_{k\in\zd}|\h{\phi^{(\beta)}}(\xi+k)|}{|\xi^\beta|}d\xi\le B'_\phi.
$$
Thus $\h\psi\in L_1$, and similarly $D^{\beta'}\h\psi\in L_1$
for all $\beta'$ such that $[\beta']\le d+1$.
 It follows that $\psi\in\mathcal{L}_p$.

 Next, we set
$$
\w\psi=\mathcal{F}^{-1}\left(\h{\w\phi}\, \gamma \right)
$$
and note that $\w\psi\in L_1$. Since $\psi\in\mathcal{L}_p$, using Lemmas~\ref{lemKK1} and~\ref{lemMZ}, we obtain
$$
\bigg\|\sum\limits_{k\in\zd}\langle T,\widetilde\psi_{0k}\rangle \psi_{0k}\bigg\|_p
\le \|\psi\|_{\mathcal{L}_p}\lll  \sum\limits_{k\in\zd} | \langle T,\widetilde\psi_{0k}\rangle|^p\rrr^{1/p}\le C_{\beta}\big\| T*\overline{\w\psi}\big\|_p= C_{\beta}\Big\|\mathcal{F}^{-1}(\h T\overline{\h{\w\psi}})\Big\|_p
$$
and
\begin{equation*}
  \begin{split}
        \Big\|\mathcal{F}^{-1}(\h T\overline{\h{\w\psi}})\Big\|_p
&=\Big\|\mathcal{F}^{-1}\lll\h T\,\overline{\h{\w\phi}}\,\gamma\rrr\Big\|_p
=\Big\|\mathcal{F}^{-1}\lll\h{D^\beta T}\,\,\overline{\h{\w\phi}}\rrr\Big\|_p\\
&=\Big\|D^\beta T*\mathcal{F}^{-1}\overline{\h{\w\phi}}\Big\|_p
\le  \Big\|\mathcal{F}^{-1}\overline{\h{\w\phi}}\Big\|_1 \|D^\beta T\|_p.
   \end{split}
\end{equation*}
Since the function $\mathcal{F}^{-1}\overline{\h{\w\phi}}$ is obviously summable on $\rd$, it follows that
\be
\label{108}
\bigg\|\sum\limits_{k\in\zd}\langle T,\widetilde\psi_{0k}\rangle \psi_{0k}\bigg\|_p\le C'_\beta \|D^\beta T\|_p.
\ee

Now let us show that
\begin{equation}\label{ad1}
  \sum\limits_{k\in\zd}\langle T,\widetilde\psi_{0k}\rangle \psi_{0k}
=\sum\limits_{k\in\zd}\langle T,\widetilde\phi_{0k}\rangle \phi^{(\beta)}_{0k}.
\end{equation}
Since the functions  $\sum_{k\in\zd}\langle T,\widetilde\psi_{0k}\rangle \psi_{0k}$ and $\sum_{k\in\zd}\langle T,\widetilde\phi_{0k}\rangle \phi^{(\beta)}_{0k}$ are locally summable, it suffices to
check that their Fourier transforms coincide almost everywhere.
Due to Lemma~1 from~\cite{KS} and Carleson's theorem, we derive
$$
\mathcal{F}\bigg(\sum_{k\in\zd}\langle  T, \w\psi_{0k}\rangle\psi_{0k}\bigg)(\xi)
=\sum_{k\in\zd}\h T(\xi+k)\overline{\h{\w\psi}(\xi+k)}\h\psi(\xi),
$$
$$
\mathcal{F}\bigg(\sum_{k\in\zd}\langle  T, \w\phi_{0k}\rangle\phi^{(\beta)}_{0k}\bigg)(\xi)
=\sum_{k\in\zd}\h T(\xi+k)\overline{\h{\w\phi}(\xi+k)}\h{\phi^{(\beta)}}(\xi).
$$
Since $ \supp \h T\subset \delta \td\subset\td$, we have for
every $\xi\in\td$ and  $l\in\zd$ that
\ban
\begin{split}
   \mathcal{F}\bigg(\sum_{k\in\zd}\langle  T, \w\psi_{0k}\rangle\psi_{0k}\bigg)(\xi+l)&=
\h T(\xi)\overline{\h{\w\psi}(\xi)}\h\psi(\xi+l)
\\
&=\h T(\xi)\overline{\h{\w\phi}(\xi)}\h{\phi^{(\beta)}}(\xi+l)=
\mathcal{F}\bigg(\sum_{k\in\zd}\langle  T, \w\phi_{0k}\rangle\phi^{(\beta)}_{0k}\bigg)(\xi+l),
\end{split}
\ean
which implies~\eqref{ad1}. Therefore, by~\eqref{108} we have
$$
I^{(\beta)}=
\bigg\|\sum\limits_{k\in\zd}\langle T,\widetilde\phi_{0k}\rangle \phi^{(\beta)}_{0k}\bigg\|_p\le C_\beta \|D^\beta T\|_p.
$$
Combining this with~\eqref{105} and~\eqref{101}, we obtain
$$
\bigg\|T-\sum\limits_{k\in\zd}\langle T,\widetilde\phi_{0k}\rangle \phi_{0k}\bigg\|_p\le C\sum_{[\beta]=s}  \|D^\beta T\|_p.
$$
To prove~\eqref{104} it remains to apply Lemmas~\ref{lem1'} and\ref{lemJ},
which completes the proof.~~$\Diamond$

\begin{theo}
\label{corMOD1'''}
Let $1\le p<\infty$, $N\ge 0$, $s\in \N$, and $\delta\in (0,1/2)$.
 Suppose that $\w\phi\in S'_{N,p}$ and $\phi\in L_p$ satisfy the following conditions:
\begin{itemize}
  \item[1)] $\supp \vp$ is compact;
  \item[2)] the Strang-Fix conditions of order $s$  hold for $\phi$;
  \item[3)] $D^{\beta}(1-\h\phi\overline{\h{\w\phi}})({\bf 0}) = 0$ for all $\beta\in\zd_+$, $[\beta]<s$;
  \item[4)] $\h{\w\vp}\in C^{s+d+1}(2\delta\td)$. 
\end{itemize}
Then, for every $f\in \mathbb{B}_{p,1}^{d/p+N}(M)$ and $j\in \z_+$, we have
\begin{equation*}
\begin{split}
\Vert f - Q_j(f,\phi,\w\phi) \Vert_p\le C \(\Omega_s(f, M^{-j})_p+m^{-j(\frac1p+\frac Nd)}\sum_{\nu=j}^\infty m^{(\frac1p+\frac Nd)\nu} E_{\delta M^\nu}(f)_p\),
\end{split}
\end{equation*}
and if, moreover, $\w\vp \in \mathcal{L}_q$, $1/p+1/q=1$,    then for every $f\in L_p$ we have
\begin{equation*}
  \Vert f - Q_j(f,\phi,\w\phi) \Vert_p\le C \Omega_s(f,M^{-j})_p,
\end{equation*}
where the constant $C$ does not depend on $f$ and $j$.
\end{theo}

{\bf Proof.}  It is sufficient to prove the theorem for   $j=0$ (see explanations in the proof of Theorem~\ref{corMOD1'+}), and  due to Lemma~\ref{thKS}, we need only to check that
\be
\label{116}
\bigg\|T-\sum\limits_{k\in\zd}\langle T,\widetilde\phi_{0k}\rangle \phi_{0k}\bigg\|_p\le
C_1\omega_s(f, 1)_p,
\ee
where   $T \in \mathcal{B}_{\delta I,p}\cap L_2$ is such that
$\Vert f-T\Vert_p\le C E_{\delta I}(f)_p$.

Let us choose a compactly supported and localy summable  function $\w\t$  such that
$$
D^{\beta}\h{\w\t}(\nul)=D^{\beta}\h{\w\phi}(\nul)\quad\text{for all}\quad \beta\in\zd_+,\,\, [\beta]<s.
$$
An appropriate function $\w\t$ can be easily constructed, e.g., as a linear combination of $B$-splines. Setting
$\w\psi=\w\phi-\w\t$, we get  $D^{\beta}\h{\w\psi}(\nul) = 0$
whenever $[\beta]<s$.

We have
\begin{equation}\label{yu1}
  \bigg\|T-\sum\limits_{k\in\zd}\langle T,\widetilde\phi_{0k}\rangle \phi_{0k}\bigg\|_p\le \bigg\|T-\sum\limits_{k\in\zd}\langle T,\w\t_{0k}\rangle \phi_{0k}\bigg\|_p+\bigg\|\sum\limits_{k\in\zd}\langle T,\w\psi_{0k}\rangle \phi_{0k}\bigg\|_p=I_1+I_2.
\end{equation}
Since both the functions $\phi$ and $\w\t$ are compactly supported, due to Theorem~\ref{Jia}  and Lemma~\ref{lemJ}, we derive
\begin{equation}\label{yu2}
  I_1\le C_2 \omega_s(T,1)_p\le C_2\(2^s \Vert f-T\Vert_p+\omega_s(f,1)_p\)\le C_3 \omega_s(f,1)_p.
\end{equation}

Let us estimate $I_2$. Since $\vp$ has a compact support and $\supp\h{T}\subset \delta\T^d$, it follows from Lemmas~\ref{lemKK1} and~\ref{lemMZ}    that
\begin{equation}\label{yu3}
  \begin{split}
     I_2^p\le C_4 \sum_{k\in \Z^d}  |\langle T, \w\psi_{0k}\rangle|^p  \le C_5 \Vert T*\w\psi\Vert_p^p.
  \end{split}
\end{equation}

Now we check that 
\be
\label{111}
\Vert T*\w\psi\Vert_p\le
C_6\sum_{[\beta]=s} \Vert D^\b T \Vert_p.
\ee
Taking again into account  that $\supp\h{T}\subset \delta\T^d$,
without loss of generality we can assume that
$\h{\w\psi}\in C^{s+d+1}(\rd)$.
By Taylor's formula with the integral reminder, we have
\begin{equation*}
  \h{\w\psi} (\xi)= \sum_{[\beta]=s} \frac{s}{\beta!}  \xi^{\beta} \int\limits_0^1 (1-t)^{s-1} D^{\beta}\h{\w\psi}(t\xi) dt.
\end{equation*}
Hence,
\begin{equation}
\label{113}
  \begin{split}
       \Vert T*\w\psi\Vert_p&=
       \left\Vert \mathcal{F}^{-1}\Big( \overline{\h{\w\psi}} \h{T}\Big)\right\Vert_{p}
\\
&\le \sum_{[\beta]=s}\frac{s}{\beta!}  \bigg(\int\limits_{\rd}\,dx\bigg| \int\limits_0^1 \, dt (1-t)^s  \int\limits_{\rd}\overline{D^\beta{\h{\w\psi}}( t\xi)}   \xi^\beta\h T(\xi) e^{2\pi i(\xi, x)}\,d\xi \bigg|^p\bigg)^{1/p}
\\
&\le\sum_{[\beta]=s}\frac{s}{\beta!}\int\limits_0^1 \, dt  \bigg(\int\limits_{\rd}\,dx \bigg|\int\limits_{\rd}\overline{D^\beta{\h{\w\psi}}( t\xi)}\h{D^\beta T}(\xi) e^{2\pi i(\xi, x)}\,d\xi\bigg|^p\bigg)^{1/p}
\\
&\le\sum_{[\beta]=s}\frac{s}{\beta!}\int\limits_0^1 \, \frac{dt}{ t^{d-d/p}} \bigg(\int\limits_{\rd}\,dx \bigg|\int\limits_{\rd}\overline{D^\beta{\h{\w\psi}}( \xi)}\h{D^\beta T}\Big(\frac \xi t\Big) e^{2\pi i(\xi, x)}\,d\xi\bigg|^p\bigg)^{1/p}
\\
&=\sum_{[\beta]=s}\frac{s}{\beta!}\int\limits_0^1 \, t^{d/p}  \bigg\Vert \mathcal{F}^{-1} \left\{\overline{D^\beta{\h{\w\psi}}}\,\h{Q_{\beta,t}} \right\}\bigg\Vert_p\, dt ,
   \end{split}
\end{equation}
where $Q_{\beta,t}(x)=D^\beta T(tx)$.
Taking into account that $D^\beta{\h{\w\psi}}\in C^{d+1}(\rd)$ whenever  $[\beta]=s$, i.e. $\mathcal{F}^{-1}\big\{D^\beta\h{\w\psi}\big\}\in L_1$, we have
\begin{equation}\label{emult}
  \begin{split}
        \bigg\Vert  \mathcal{F}^{-1}\left\{\overline{D^\beta{\h{\w\psi}}}\,\h{Q_{\beta,t}} \right\}\bigg\Vert_p
	&=\|\mathcal{F}^{-1}\big\{D^\beta\h{\w\psi}\big\} * Q_{\beta,t}\|_p
	\\
	&\le\|\mathcal{F}^{-1}\big\{D^\beta\h{\w\psi}\big\}\|_{1}  t^{-d/p}\|D^\beta T\|_p=C_8  t^{-d/p}\|D^\beta T\|_p.
   \end{split}
\end{equation}
Combining these estimates with~\eqref{113}, we obtain~\eqref{111}.  Next, similarly as in~\eqref{yu2}, using~\eqref{111}, \eqref{yu3}, and
Lemmas~\ref{lem1'} and~\ref{lemJ}, we get
\begin{equation}\label{fine}
  I_2 \le C_9 \omega_s(f,1)_p.
\end{equation}

Finally, combining~\eqref{yu1}, \eqref{yu2}, and~\eqref{fine}, we complete the proof.~~$\Diamond$

\begin{rem}
The conditions on smoothness of $\h{\w\vp}$ in Theorem~\ref{corMOD1'''} can be relaxed and given in other terms using different sufficient conditions for Fourier multiplier in $L_p$ spaces, see e.g.,~\cite{LST} and~\cite{Kol}. For this, one needs to estimate the corresponding multiplier norm of the function $D^\beta\h{\w\psi}$ instead of to show that $\mathcal{F}^{-1}\big\{D^\beta\h{\w\psi}\big\}\in L_1$ in inequality~\eqref{emult}. A similar conclusion is valid also for the smoothness conditions on $\h\vp$ and $\h{\w\vp}$ in Theorem~\ref{corMOD1'+}.
\end{rem}

\begin{rem}
It is not difficult to verify that Theorems~\ref{corMOD1'+}  and~\ref{corMOD1'''} are valid also in the spaces $L_\infty$ if we additionally suppose  that  $f\in L_2\cap L_\infty$ and replace the best approximation $E_{\delta M^\mu}(f)_p$ with $E_{\delta M^\mu}^*(f)_p$. See also~\cite[Theorem 17]{KSa} for the corresponding analogue of Lemma~\ref{lemJ}.
\end{rem}

Let us compare  Theorems~\ref{corMOD1'+}  and~\ref{corMOD1'''}  with Theorem~\ref{theoQj} in the case of the isotropic matrix $M$, $\w\phi\in S'_{N,p}$, and  $p\ge2$.  First let us show that the class of
functions $f$ considered in Theorem~\ref{theoQj} is smaller than  $\mathbb{B}_{p,1}^{d/p+N}(M)$. Let
$\l$ be an eigenvalue of $M$, $f\in L_p$, $\h f\in L_q$, and
 $\h f(\xi)=O(|\xi|^{-N-d-\varepsilon})$ as $|\xi|\to\infty$, $\varepsilon>0$.
  Setting  $V_\mu f=\mathcal{F}^{-1}v *{f}$,
where $v\in C^\infty(\R)$, $v(\xi)\le1$,  $v(\xi)=1$ if $|M^{-\mu}\xi|\le 1/2$ and $v(\xi)=0$ if $|M^{-\mu}\xi|\ge 1$, using
Pitt's inequality (see, e.g., \cite[inequality (1.1)]{GT})
$$
\Vert f\Vert_{p}\le C(p)\bigg(\,\int\limits_{\R^d} |\xi|^{d(p-2)}|\h f(\xi)|^p d\xi\bigg)^\frac1p,\quad 2<p<\infty,
$$
and taking into account~\eqref{10}, we obtain
\begin{equation*}
  \begin{split}
      E_{M^\mu}(f)_p&\le \Vert f-V_\mu(f)\Vert_p
			\le C(p)\(\,\int\limits_{|M^{-\mu}\xi|\ge 1/2}|\xi|^{d(p-2)} |(1-v(\xi))\h f(\xi)|^p d\xi\)^\frac1p
			\\
      &\le 2C(p)\(\,\int\limits_{|\xi|\ge 1/2\|M^{-\mu}\|}|\xi|^{d(p-2)} |\h f(\xi)|^p d\xi\)^\frac1p
      =\mathcal{O}\(|\lambda|^{-\mu(d/p+N+\varepsilon)}\),\quad \mu\to\infty.
   \end{split}
\end{equation*}
It follows that  $f\in \mathbb{B}_{p,1}^{d/p+N}(M)$.

Next let us compare the error estimates. Using~\eqref{10}, we have
$$
\Omega_s(f, M^{-j})\le C_1\omega_s(f, |\lambda|^{-j}),
$$
where $C_1$ does not depend on $f$ and $j$, and
$$
E_{|\lambda|^\mu I}(f)\le C_2  E_{\sigma M^\mu}(f)_p\le C_3 E_{M^{\mu-\mu_0}}(f)_p,
$$
where  $C_2$, $C_3$, $\sigma$, and $\mu_0$ do not depend on $f$ and $j$. Hence
$$
E_{|\lambda|^\mu I}(f)_p =\mathcal{O}\(|\lambda|^{-\mu(d/p+N+\varepsilon)}\),\quad \mu\to\infty.
 $$
 Using these relations and the following inverse approximation inequality (see~\cite{DDT}):
\begin{equation}\label{eqMarch}
  \omega_s(f,2^{-j})_p\le C2^{-s j}\(\sum_{k=0}^j 2^{2s k} E_{{2^k}I}(f)_p^2\)^\frac12,
\end{equation}
we can easily see that  inequality~\eqref{110} provides
the following approximation order
	 $$
\Vert f - Q_j(f,\phi,\w\phi) \Vert_p\le
	 \begin{cases}
	 C |\lambda|^{-j(N+\frac dp + \varepsilon)}  &\mbox{if }
	s> N+\frac dp + \varepsilon\\
	  C (j+1)^{1/2} |\lambda|^{-js} &\mbox{if }
	 s= N+\frac dp + \varepsilon \\
	C|\lambda|^{-js}
	 &\mbox{if }
	 s< N+\frac dp + \varepsilon
	\end{cases},
$$
which is better than one given in Theorem~\ref{theoQj}  in the case $s= N+\frac dp + \varepsilon $,
and the same in the other cases.

On the other hand, there exist functions in $\mathbb{B}_{p,1}^{d/p+N}(M)$ which do not satisfy assumptions of Theorem~\ref{theoQj}.
Indeed, let  $b\le N+d$,  $\varepsilon >0$,
$$
a>\max\left\{1,\   (d+N+\varepsilon-b)\lll\frac12-\frac1p\rrr^{-1}\right\},
$$
$$
\h f(\xi) = \kappa(|\xi|)\frac{e^{i|\xi|^a}}{|\xi|^{b}},
$$
where $\kappa\in C^{\infty}(\Bbb R)$, $ \kappa(u)=0$ for $u=[0,1]$,
$ \kappa(u)=1$ for $u\ge2$.
Obviously, the decay of $\h f$ is not enough for Theorem~\ref{theoQj}.
 Let us verify that $f:=\mathcal{F}^{-1}(\h f)$ is in $L_p$ and
$E_{M^\nu}(f)_p=\mathcal{O}(|\lambda|^{-\gamma \nu})$,
where  $\gamma=\frac dp+N+\varepsilon$, i.e., $f\in \mathbb{B}_{p,1}^{d/p+N}(M)$.
Setting $g(\xi):=|\xi|^\gamma\h f(\xi)$ and using Proposition~5.1 from~\cite{Mi}) we conclude that the functions
$\mathcal{F}^{-1} \h f$ and  $\mathcal{F}^{-1} g$ are smooth throughout $\R^d$ and
$$
\mathcal{F}^{-1}(\h f)(x)=d_0 |x|^{-\frac{b-d+da/2}{a-1}}\exp\(id_1 |x|^{\frac{a}{a-1}}\)+o\(|x|^{-\frac{b-d+da/2}{a-1}}\)\quad\text{as}\quad |x|\to \infty,
$$
$$
\mathcal{F}^{-1}(g)(x)=d_0 |x|^{-\frac{b-\gamma-d+da/2}{a-1}}
\exp\(id_1 |x|^{\frac{a}{a-1}}\)+o\(|x|^{-\frac{b-\gamma-d+da/2}{a-1}}\)\quad\text{as}\quad |x|\to \infty.
$$
It follows that $f\in L_p$ and $\mathcal{F}^{-1} g\in L_p$, which yields that  $f\in \dot H_p^\gamma$ (the homogeneous Sobolev spaces). Then, using the embedding
$\dot H_p^\gamma \subset \dot B_{p,\infty}^\gamma$ (the homogeneous Besov space, see, e.g.,~\cite[Theorem~6.3.1]{BL}) and  Jackson's theorem, we obtain
$E_{M^\nu}(f)_p\le C\omega_{[\gamma]+1}(f,|\lambda|^j)_p\le C\Vert f\Vert_{\dot B_{p,\infty}^\gamma}|\lambda|^{-\gamma j} $.


Finally we consider an anisotropic case. Let $d=2$, $M={\rm diag} \(m_1,m_2\)$, and $\frac{\partial^{r_1}}{\partial x_1^{r_1}} f, \frac{\partial^{r_2}}{\partial x_2^{r_2}} f \in L_p$, where $r_1>(\frac1p+\frac Nd)(1+\frac{\ln m_2}{\ln m_1})$ and $r_2>(\frac1p+\frac Nd)(1+\frac{\ln m_1}{\ln m_2})$. Then, taking into account that for $1<p<\infty$ and $s>\max(r_1,r_2)$, see, e.g.,~\cite{timanIz}
$$
E_{M^\nu}(f)_p\le C \Omega_s(f,M^{-\nu})_p\le C\(\omega_{s}^{(1)}(f,m_1^{-j})_p+\omega_{s}^{(2)}(f,m_2^{-j})_p\)=\mathcal{O}\({m_1^{-r_1 j}}+{m_2^{-r_2 j}}\),
$$
where $\omega_{s}^{(\ell)}(f,h)_p$ is the partial modulus of smoothness with respect to $\ell$-th variable,
we have that
$f\in \mathbb{B}_{p,1}^{d/p+N}(M)$, and by Theorem~\ref{corMOD1'+} with $s>\max(r_1,r_2)$, we get
$$
\Vert f-Q_j(f,\vp,\w\vp)\Vert_p=\mathcal{O}\({m_1^{-r_1 j}}+{m_2^{-r_2 j}}\).
$$
That is, the order of approximation by $Q_j(f,\vp,\w\vp)$ essentially depends on anisotropic nature of the function $f$ and the matrix $M$.

\section{Conclusions}
Error estimates in $L_p$-norm  are obtained for a large class of
sampling-type quasi-projection operators $Q_j(f,\phi, \w\phi)$ including
the classical case, where  $\w\phi$ is the Dirac delta-function.
Theorems~\ref{corMOD1'+} and~\ref{corMOD1'''} provide  essential improvements of the recent results in~\cite{KS}. First, the estimates are given for a wider class of approximated functions, namely for functions from anisotropic Besov spaces.  Second, only the case $p\ge2$  was considered in~\cite{KS}, while  $1\le p<\infty$ in Theorems~\ref{corMOD1'+} and~\ref{corMOD1'''}.
Third, the estimates are given in terms of the moduli of smoothness and best approximations,
while the results in~\cite{KS} provide only approximation order.

\bigskip

\noindent {\bf ACKNOWLEDGMENTS}

The first author was partially supported by DFG project KO 5804/1-1 (Theorem~\ref{corMOD1'''} and the corresponding auxiliary results belong to this author).
The second author was supported by the Russian Science Foundation under grant No. 18-11-00055 (Theorem~\ref{corMOD1'+} and the corresponding auxiliary results belong to this  author).


\begin{thebibliography}{}


\bibitem{BL} {\sc J. Bergh, J. L\"ofstr\"om}, Interpolation spaces. An introduction, Springer, 1976.



\bibitem{Butz4}
        {\sc C. Bardaro, P.\,L. Butzer, R.\,L. Stens, G. Vinti},
        Approximation error of the Whittaker cardinal series in terms of an
       averaged modulus of smoothness covering discontinuous signals,
       {\it Math. Anal. Appl.} {\bf 316} (2006), no.~1, 269--306.
			
			\bibitem{Butz6}
        {\sc C. Bardaro, P.\,L. Butzer, R.\,L. Stens, G. Vinti},
       Prediction by samples from the past with error estimates covering
       discontinuous signals,
        {\it IEEE Trans. Inform. Theory} {\bf 56} (2010), no.~1, 614--633.



\bibitem{BDR} {\sc C. de Boor, R. DeVore, A. Ron}, Approximation from shift-invariant subspaces of $L_2(\R^d)$, {\it Trans. Amer. Math. Soc.} \textbf{341} (1994), no.~2, 787--806.

\bibitem{DB-DV1} {\sc C. de Boor, R. DeVore, A. Ron},
Approximation orders of FSI spaces in $L_2(\rd)$,
{\it Constr. Approx.} \textbf{14} (1998), no.~3, 411--427.


\bibitem{Brown}
        {\sc J.\,L. Brown, Jr.},
        On the error in reconstructing a non-bandlimited function by means of
        the bandpass sampling theorem,
        {\it J. Math. Anal. Appl.} {\bf 18} (1967), 75--84.


\bibitem{BD} {\sc M. D. Buhmann, F. Dai}, Pointwise approximation with quasi-interpolation by radial basis functions, {\it J. Approx. Theory} \textbf{192} (2015), 156--192.




\bibitem{Butz5}
       {\sc P.\,L. Butzer, J.\,R. Higgins, R.\,L. Stens},
        Classical and approximate sampling theorems: studies in the
       $L_p(\mathbb{R})$ and the uniform norm,
       {\it J. Approx. Theory} {\bf 137} (2005), no. 2, 250--263.


\bibitem{Butz7}
        {\sc P.\,L. Butzer and R.\,L. Stens},
        Reconstruction of signals in $L_p(\mathbb{R})$-space by generalized
        sampling series based on linear combinations of B-splines,
        {\it Integral Transforms Spec. Funct.} {\bf 19} (2008), no. 1, 35--58.



\bibitem{CV0} {\sc D. Costarelli, G. Vinti}, Approximation by nonlinear multivariate sampling Kantorovich type operators and
applications to image processing, \emph{Numer. Funct. Anal. Optim.} \textbf{34} (2013) 819--844.


\bibitem{CV2} {\sc D. Costarelli, G. Vinti}, Rate of approximation for multivariate sampling Kantorovich operators
on some functions spaces, {\it J. Int.
Eq. Appl.} \textbf{26} (2014), no.~4,  455--481.

\bibitem{DDT} {\sc F. Dai, Z. Ditzian, S. Tikhonov}, Sharp Jackson inequalities, {\it J. Approx. Theory} {\bf 151} (2008), 86--112.




\bibitem{GT}
{\sc L. De Carli, D. Gorbachev,  S.  Tikhonov}, Pitt and Boas inequalities for Fourier and Hankel transforms. {\it J. Math. Anal. Appl.} {\bf 408} (2013), no. 2, 762--774.

 \bibitem{v58}
{\sc R.-Q. Jia}, Refinable shift-invariant spaces: from splines to wavelets, in: C.K. Chui, L.L. Schumaker (Eds.), Approximation Theory VIII, vol. 2 (College Station, TX, 1995), Ser. Approx. Decompos., vol. 6, World Scientific Publishing, River Edge, NJ, 1995, pp. 179--208.


\bibitem{Jia1} {\sc R.-Q. Jia}, Convergence rates of cascade algorithms, {\it Proc. Amer. Math. Soc.} \textbf{131} (2003), 1739--1749.

\bibitem{Jia2}
{\sc R.-Q. Jia}, Approximation by quasi-projection operators in Besov spaces,
{\it J. Approx. Theory} \textbf{162} (2010), no.~1, 186--200.


\bibitem{Kol}
{\sc Yu. S. Kolomoitsev}, Multiplicative sufficient conditions for Fourier multipliers, {\it Izv. Math.} \textbf{78} (2014), no.~2,   354--374.


\bibitem{KKS} {\sc  Yu. Kolomoitsev, A. Krivoshein, M. Skopina}, Differential and falsified sampling expansions, \emph{J. Fourier Anal. Appl.} \textbf{24} (2018), no.~5, 1276--1305.




\bibitem{KS3} {\sc Yu. Kolomoitsev, M. Skopina}, Approximation by multivariate Kantorovich-Kotelnikov operators, \emph{J. Math. Anal. Appl.} \textbf{456} (2017), no.~1, 195--213.


\bibitem{KSa} {\sc Yu. Kolomoitsev, M. Skopina}, Quasi-projection operators in the weighted $L_p$ spaces, arXiv:1805.10536v1.



\bibitem{KS}  {\sc A. Krivoshein, M. Skopina}, Multivariate sampling-type approximation, \emph{Anal. Appl.} \textbf{15} (2017), no.~4, 521--542.




\bibitem{LST}   {\sc E. Liflyand, S. Samko, R. Trigub}, The Wiener algebra of absolutely convergent Fourier
integrals: an overview, \emph{Anal. Math. Phys.} \textbf{2} (2012), 1--68.


\bibitem{Mi} {\sc A. Miyachi}, On some singular Fourier multipliers, \emph{J. Fac. Sci. Univ. Tokyo Sect. IA Math.} \textbf{28} (1981), 267--315.


\bibitem{NU} {\sc H. Q. Nguyen, M. Unser}, A sampling theory for non-decaying signals, \emph{Appl. Comput. Harmon. Anal.} \textbf{43} (2017), no.~1, 76--93.





\bibitem{Nik} {\sc S. M. Nikol'skii},
The Approximation of Functions of Several Variables and the Imbedding Theorems, 2nd ed.
{Moscow}: Nauka, 1977 (Russian). -- English transl. of 1st. ed.: John Wiley $\&$ Sons, New-York, 1978.



\bibitem{OT15} {\sc O. Orlova, G. Tamberg},  On approximation properties of generalized
Kantorovich-type sampling operators, \emph{J. Approx. Theory} {\bf 201} (2016), 73--86.


\bibitem{SS} {\sc H.-J. Schmeisser, W. Sickel}, Sampling theory and function spaces. In : Applied Mathematics Reviews, Vol. 1, 205-284, World Scientific, 2000.



				
\bibitem{Si2} {\sc W. Sickel}, Spline representations of functions in Besov-Triebel-Lizorkin spaces on $\mathbb{R}^n$,
				{\it Forum Math.} \textbf{2} (1990), no.~5, 451--475.


				

\bibitem{Sk1}
         {\sc M. Skopina},
        Band-limited scaling and wavelet expansions,
        {\it Appl. Comput. Harmon. Anal.} {\bf 36} (2014), 143--157.




\bibitem{Timan} {\sc A. F. Timan}, Theory of Approximation of Functions of a Real Variable, Pergamon Press, Oxford, London, New York, Paris, 1963.

\bibitem{timanIz} {\sc M. F. Timan}, The difference properties of functions of several variables, (Russian) Izv. Akad. Nauk SSSR Ser. Mat. \textbf{33} (1969) 667--676.


\bibitem{TB} {\sc R. M. Trigub, E. S. Belinsky}, Fourier Analysis and Appoximation of Functions. Kluwer. 2004.



\bibitem{Unser}
{\sc M. Unser}, Sampling - 50 years after Shannon, \emph{Proceedings of the IEEE} \textbf{88} (2000), 569--587.



\bibitem{VZ2} {\sc G. Vinti, L. Zampogni}, Approximation results for a general class of Kantorovich type operators, \emph{Adv. Nonlinear
Stud.} \textbf{14 }(2014), no.~4, 991--1011.



\bibitem{Wil} {\sc  G. Wilmes},  On Riesz-type inequalities and $K$-functionals related to Riesz potentials in $\R^N$,
{\it Numer. Funct. Anal. Optim.} {\bf 1} (1979), no.~1, 57--77.



\end{thebibliography}
\end{document}